\documentclass[letterpaper, 10 pt, conference]{ieeeconf}  % Comment this line out if you need a4paper

\IEEEoverridecommandlockouts                              % This command is only needed if 
\usepackage{cite}
\usepackage{amsmath,amssymb,amsfonts}
\usepackage{graphicx}
\usepackage{textcomp}

\usepackage{graphics} % for pdf, bitmapped graphics files
\usepackage{epsfig} % for postscript graphics files
\usepackage{mathptmx} % assumes new font selection scheme installed
\usepackage{times} % assumes new font selection scheme installed
\usepackage{amsbsy}
\usepackage{dsfont}
\usepackage{verbatim}
\usepackage{algorithm}
\usepackage[colorlinks=true,urlcolor=blue]{hyperref}
\usepackage{algpseudocode}
\usepackage{url}
\usepackage{xcolor}
\graphicspath{{./images/}}

\newcommand{\eps}{\varepsilon}

\DeclareMathOperator*{\argmin}{arg\,min}

\newcommand{\R}{\mathbb{R}}
\newcommand{\B}[1]{\mathbf{#1}}

\newcommand{\abs}[1]{\left \vert #1 \right \vert}
\algrenewcommand\algorithmicindent{1em}

\begin{document}

\title{
Collisionless Multi-Agent Path Planning in the Hamilton-Jacobi Formulation
}
\author{Christian Parkinson, Adan Baca and Huy Nguyen
\thanks{*Early stages of this work were supported by the NSF through the Research Training Group on Applied Mathematics and Statistics for Data Driven Discovery at the University of Arizona (DMS-1937229).}
\thanks{Christian Parkinson is an assistant professor in the Department of Mathematics and Department of Computational Mathematics, Science and Engineering, Michigan State University, East Lansing, MI, USA, 48824
        {\tt\small chparkin@msu.edu}}%
\thanks{Adan Baca is a student in the Department of Mathematics and Department of Computer Science, University of Arizona, Tucson, AZ, USA, 85721
        {\tt\small adanbaca@arizona.edu}}%
\thanks{Huy Nguyen is a student in the Department of Mathematics, Michigan State University, East Lansing, MI, USA, 48824 
        {\tt\small nguy1270@msu.edu}}%
}

\maketitle

\begin{abstract}
We present a method for collisionless multi-agent path planning using the Hamilton-Jacobi-Bellman equation. Because the method is rooted in optimal control theory and partial differential equations, it avoids the need for hierarchical planners and is black-box free. Our model can account for heterogeneous agents and realistic, high-dimensional dynamics. We develop a grid-free numerical method based on a variational formulation of the solution of the Hamilton-Jacobi-Bellman equation which can resolve optimal trajectories even in high-dimensional problems, and include some practical implementation notes. In particular, we resolve the solution using a primal-dual hybrid gradient optimization scheme. We demonstrate the method's efficacy on path planning problems involving simple cars and quadcopter drones.
\end{abstract}

%\begin{IEEEkeywords}
%Optimal Path-Planning, Multi-Agent, Coordination, Dynamic Programming, Hamilton-Jacobi Equation, Hopf-Lax Formula 
%\end{IEEEkeywords}

\thispagestyle{empty}

%%%%%%%%%%%%%%%%%%%%%%%%%%%%%%%%%%%%%%%%%%%%%%%%%%%%%%%%%%%%%%%%%%%%%%%%%%%%%%%%
\section{INTRODUCTION}

Multi-agent path planning in crowded environments is a problem of keen interest to the automation and control community with applications like building evacuation \cite{Evac}, reduction of cluttered air corridors \cite{AirTravel}, and optimal terrain coverage \cite{Coverage}. In this manuscript, we present a partial differential equation (PDE) based method for time-optimal collisionless multi-agent path planning based on optimal control theory, dynamic programming, and a Hamilton-Jacobi-Bellman (HJB) equation. 

With advances in computing power over the past years, there are many modern multi-agent planning methods based on hierarchical algorithms and deep reinforcement learning \cite{DL1,DL2,DL3}. Deep learning methods are applicable to high-dimensional problems and exceedingly flexible with respect to complex dynamics and constraints. However, as is now well documented, explainability and robustness of deep learning methods are sometimes questionable, and while addressing these concerns is an active area of research across several fields \cite{Robust1,Explain1,Robust2,Explain2}, PDE-based planning methods provide a fully explainable and robust alternative. 

Path planning methods which employ a HJB formulation have been successfully used  in the context of underwater navigation with uncertain currents \cite{Lolla}, control theoretic models of environmental crime \cite{Arnold, Cartee2, Chen}, human-navigation in mountainous terrain \cite{Parkinson,Parkinson2}, surveillance-evasion games \cite{Cartee1}, animal foraging \cite{GeeFear}, and single-agent planning for simple autonomous vehicles \cite{TakeiTsai1,TakeiTsai2,ParkinsonCar1,ParkinsonCar2,Gee}, to list a few. All the work contained in these references takes place, by necessity, in low-dimensional state spaces, where the relevant PDE can be solved using classical grid based methods such as fast sweeping and fast marching type schemes \cite{FS1,FS4,Sethian,SethVlad1}. Historically, this has been the primary drawback of PDE-based path planning methods: the classical numerical schemes suffer from the curse of dimensionality, limiting their utility to low dimensional problems.  

There has been much recent work to ameliorate this drawback through the design of scalable numerical methods for Hamilton-Jacobi equations (see \cite{HighDimHJ} for a summary and review of many such methods). One particular class of such methods looks for characteristics-based solutions via Hopf-Lax type representations of the solution of HJB equations \cite{Lin,ParkinsonBoyle,ParkinsonPolage,ParkinsonBaca}. In this vein, we present a method for computing optimal trajectories for minimal-time collisionless multi-agent path planning in obstacle ridden domains by computing characteristic curves of an associated HJB equation. Our method is amenable to heterogeneous agents and efficient enough to be employed in semi-real-time, accounting for a scenarios in which new obstacles are discovered as the agents traverse the domain, and for reasonably complex dynamics like those of a quadcopter drone.

\section{OPTIMAL CONTROL AND THE HAMILTON-JACOBI-BELLMAN EQUATION}

In this section, we describe the optimal control problem that we will solve and formally derive the HJB equation. Before getting to this, we establish some notational conventions. We consider trajectory planning for $I$ agents. The subscript $i$ will always refer to the $i^{\text{th}}$ agent. Throughout, trajectories in state or control space will be represented by boldface letters, while points in state or control space will be in normal font. For brevity, we use brackets rather than tuples to represent the position of the ensemble in the holistic state space. Thus, $\{x_i\}$ is taken to mean $(x_1,x_2,\ldots,x_I)$. 

Suppose we have agents with trajectories represented by $\{\B x_i(\cdot)\}_{i=1}^I$. In the absence of obstacles or other agents, agent $i$ obeys dynamics \begin{equation} \label{eq:dynamics} 
\dot{\B x}_i = f_i(\B x_i, \B a_i, t). 
\end{equation} Allowing for heterogeneous agents, we assume that dynamics functions $f_i:\R^{n_i}\times \R^{m_i}\times [0,T] \to \R^{n_i}$ and initial conditions $x_{i,0} \in \R^{n_i}$ are given, and control trajectories are maps $\B a_i: [0,T] \to \mathbb A_i\subset \R^{m_i}$, where the state dimensions $n_i$ and control dimensions $m_i$ for different agents are possibly different.  Thus the holistic dynamics take place in the $\sum_i n_i$ dimensional space $\Pi_i \R^{n_i}$ and the holistic controls are taken from the $\sum_i m_i$ dimensional space $\Pi_i \mathbb A_i$. We assume collaboration and perfect sharing of information between the agents so that this can be considered a centralized control problem. Here $T > 0$ is the horizon time, which is somewhat arbitrary, as we address shortly. We further suppose that we are given a desired final configuration $\{x_{i,f}\} \in \Pi_i \R^{n_i}$ such that the goal of the controller is to navigate the agents to this final position in minimal time. To model this mathematically, for fixed $\{x_i\} \in \Pi_i \R^{n_i}$ and $t \in [0,T]$, we define the remaining cost functional \begin{equation} \label{eq:cost}
\mathcal C_{\{x_i\},t}[\{\B x_i\}, \{\B a_i\}] = g(\{\B x(T)\}) + \int^T_t \chi_{f}(\{\B x_i(s)\})ds.
\end{equation} where we have restricted ourselves to trajectories $\{\B x_i(\cdot)\}$ such that $\B x_i(t) = x_i$, and control values $\{\B a_i(s)\}$ for $s \in [t,T]$. Here $g(\{x_i\}) = 0$ if $\{x_i\}= \{x_{i,f}\}$ and $g(\{x_i\}) = +\infty$ otherwise. Likewise, $\chi_f(\{x_i\}) = 0$ if $\{x_i\}= \{x_{i,f}\}$ and $\chi(\{x_i\}) = 1$ otherwise. Hence to minimize this cost functional, the controller must steer the agents to their desired final configuration by the final time to avoid an infinite cost. Among such paths, the running cost, represented by the integral in \eqref{eq:cost}, will count marginal cost of $1$ per unit time until it ``shuts off" when all agents reach their final configuration. Thus minimizing this cost functional entails steering the agents to their final configurations as quickly as possible. 

The value function is then defined by \begin{equation} \label{eq:value} u(\{x_i\},t) = \inf_{\{\B a_i(s)\}, t \le s \le T} \mathcal C_{\{x_i\},t}[\{\B x_i(\cdot)\}, \{\B a_i(\cdot)\}]. \end{equation} In words, $u(\{x_i\},t)$ is the minimal remaining cost achievable for a trajectory that is at $\{x_i\}$ at time $t$.  Assuming that agents are at $\{x_{i,0}\}$ at $t = 0$, the goal is then to resolve $u(\{x_{i,0}\},0)$ as well as the controls $\{\B a_i(\cdot)\}$ which achieve this value in \eqref{eq:value}. A classical result \cite{Bardi,Tran} states that the value function is the viscosity solution of the HJB equation \begin{equation}
\label{eq:HJB0}
\begin{aligned} &u_t + \inf_{\{a_i\}}\left\{ \left( \sum^I_{i=1} \langle f_i(x_i,a_i,t),\nabla_i u\rangle\right)+ \chi_f(\{x_i\})\right\}=0, \\
&u(\{x_i\},T) = g(\{x_i\}), \end{aligned}
\end{equation} where $\nabla_i$ denotes the gradient of $u$ with respect to the variables corresponding to agent $i$. 

We now make several small modifications to this equation for computational convenience. Recall, the dynamics \eqref{eq:dynamics} were stated disregarding collisions and obstacles. A typical manner of incorporating these is to restrict decompose the state space $\Pi_i \R^{n_i}$ into legal configurations---those for which no agents are colliding with each other or obstacles---and illegal configurations---those for which at least two agents are colliding or at least one agent is colliding with an obstacle (see e.g. \cite{ParkinsonCar2}). This leads to a constrained state space or, equivalently, some spatial ``boundary" conditions that $u(\{x_i\},t)$ must satisfy. Numerically, we will solve the HJB equation using a variational formulation, and these conditions become unwieldy constraints in the optimization routines. Accordingly, we would like to avoid explicitly enforcing any spatial conditions. These conditions can instead be implicitly enforced using a strategy similar to that in \cite{TakeiTsai2,ParkinsonBoyle,ParkinsonPolage}. The basic idea is to set the velocity of any agents that have collided with each other or obstacles to zero. Having done this, agents which take illegal paths will get stuck, preventing them from reaching their final destination, and incurring an infinite cost from the terminal cost function.  

To accomplish this, we define pairwise collision function $c_{k\ell}(x_k,x_\ell)$ which take the values $0$ when agents $k$ and $\ell$ are colliding and $1$ otherwise, and obstacle collision functions $O_i(x_i,t)$ which take the value $0$ when agent $i$ is colliding with an obstacle at time $t$ and $1$ otherwise. The holistic collision function is then given by $C(\{x_i\}) = \Pi_{k=1}^I \Pi_{\ell=k+1}^I c_{k\ell}(x_k,x_\ell)$. We then multiply the dynamics functions with these collision functions. For one last computational trick, we note that assuming $\B x_i(t) = x_{i,f}$, the optimal strategy $\B a_i(t)$ will certainly stop the motion for agent $i$. We can avoid explicitly accounting for this spatial condition, but instead multiplying the dynamics by $\chi_f(\{x_i\})$ which will force the agents to stop when all have reached their final points. With these considerations, we actually use the dynamics $$\dot{\B x}_i =  \chi_f(\{\B x_k\})C(\{\B x_k\}) O_i(\B x_i,t) f_i (\B x_i, \B a_i,t).$$  Lastly, because initial value problems are more comfortable for those working in PDE, and because it allows us to eschew the horizon time $T$, we make the time reversing substitution $t\mapsto T-t$. Incorporating all of these modifications, \eqref{eq:HJB0} becomes \begin{equation}\label{eq:HJB}\begin{aligned} &u_t + \mathcal H\Big(\{x_i\},\{\nabla_iu\},t\Big) = 0, \\ &u(\{x_i\},0) = g(\{x_i\}), \end{aligned} \end{equation} where, using $p_i$ as a proxy for $\nabla_i u$ and suppressing arguments for brevity, the holistic Hamiltonian is defined \begin{equation} \label{eq:Ham} \mathcal H\Big(\{x_i\},\{p_i\},t\Big) = \chi_f\left( C \left(\sum^I_{i=1} O_i H_i\right) - 1 \right) \end{equation} with $\chi_f(\{x_i\}), C(\{x_i\})$ and $O_i(x_i,t)$ as above, and individual Hamiltonians defined \begin{equation}\label{eq:Hi}H_i(x_i,p_i,t) = \sup_{a_i} \langle -f_i(x_i,a_i,t),p_i\rangle.\end{equation} Assuming the solution of \eqref{eq:HJB} is known, the optimal feedback control action $\{a_i\}$ for agents in configuration $\{x_i\}$ is the action that achieves the supremum for each of these individual Hamiltonians.

Having made the time-reversing substitution, $T$ no longer explicitly appears in the HJB equation, and $t$ now represents \emph{remaining travel time}. Note that at $t = 0$, the solution $u(\{x_i\},0) = g(\{x_i\})$ is infinite except at $\{x_{i,f}\}$, denoting that if there is zero travel time remaining, the controller is assessed infinite cost unless the agents reached their final points. At later times, $u(\{x_i\},t) < +\infty$ if and only if there is an admissable control which steers the agents from $\{x_i\}$ to $\{x_{i,f}\}$ in time $\le t$. We make two basic assumptions about the dynamics: (1) the system has small-time local controllability \cite{STLC}, meaning in essence, that if $\|\{x_i\} -\{y_i\}\| \le \eps$, then the agents can reach $\{y_i\}$ from $\{x_i\}$ in $O(\eps)$ time, and (2) from any configuration $\{x_i\}$, the agents can reach the final configuration $\{x_{i,f}\}$ in finite time. Given these, for any compact set $K\subset \Pi_i \R^{n_i}$, there is a maximal time $T^*$ such that for any $\{x_i\} \in K$, the agents can reach $\{x_{i,f}\}$ in time $\le T^*$. In what follows, we develop a numerical scheme to approximately solve \eqref{eq:HJB} at individual points $(\{x_i\},t)$. Due to the previous observation, we are justified in restricting $\{x_i\}$ to a compact set and fixing $t$ large enough. This essentially removes the time horizon $T$ from the problem, assuming it is large enough, while still allowing for time-dependent dynamics if desired. 

In the ensuing section, we describe a grid free numerical method for approximating the solution of \eqref{eq:HJB} while also explicitly resolving optimal trajectories. 

\section{NUMERICAL METHODS}

In this section, we present a method for computationally approximating the solution \eqref{eq:HJB} at individual points $(\{x_i\},t)$. The general strategy is to write $u(\{x_i\},t)$ as the solution of a saddle point problem for an appropriately defined Lagrangian, in a manner similar to the classical Hopf-Lax formula \cite[Thm. 2.23]{Tran}. Assuming the Hamiltonian is space- and time-independent (and under further mild conditions on $\mathcal H$ and $g$), the classical Hopf-Lax formula gives the solution of $u_t + \mathcal H(\nabla u) = 0, \,\, u(x,0) = g(x)$ via the saddle point problem \begin{equation}\label{eq:classicHL} u(x,t) = \inf_{\overline x} \sup_{\overline p} \Big\{g(\overline x) + \langle x-\overline x, \overline p\rangle - tH(\overline p)\Big\}.  \end{equation} However, for most interesting optimal control problems (as with ours), the Hamiltonian will depend on space and/or time. For such problems, the classical Hopf-Lax formula no longer holds. Using similar methodology, the authors of \cite{Lin,ParkinsonBoyle,ParkinsonPolage,ParkinsonBaca} have proposed a conjectural discrete Hopf-Lax formula which applies for space- and time-dependent Hamiltonians. We describe this method, and its adaptation it to our particular problem. For more details, see \cite{Lin,ParkinsonBoyle}. 

Fixing a point $(\{x_i\},t)$, we let $0 = t_0 < t_1 < \ldots < t_J$ be a uniform discretization of the time interval $[0,t]$ with step size $\Delta t = t/M$. Denote by $x_{i,j}$ our discrete approximation to $\B x_i(t_j)$. We introduce discrete co-state vectors $p_{i,j}$ corresponding to agent $i$ and time step $j$. These can be realized as an approximation to the optimal co-state dynamics ensured by the Pontryagin maximum principle \cite[\S1.6]{Bardi} or equivalently as approximations to $\nabla_i u(\{\B x_i(t_j)\},t_j)$ whenever the solution of \eqref{eq:HJB} is smooth. As described in \cite{Lin,ParkinsonBoyle}, mathematically they enter into the problem as Lagrange multipliers which enforce a discrete version of the desired dynamics.With this notation, the authors of \cite{Lin,ParkinsonBoyle,ParkinsonPolage,ParkinsonBaca} conjecture (and provide solid empirical evidence) that the solution of \eqref{eq:HJB} can be approximated using \begin{equation} \label{eq:HL} \begin{aligned}
u(\{x_i\},t) \approx \inf_{\{\overline x_{i,j}\}} \sup_{\{\overline p_{i,j}\}} \Big\{ g&(\{\overline x_{i,0}\}) + \sum_{i,j}  \left\langle \overline x_{i,j+1} - \overline x_{i,j}, \overline p_{i,j} \right \rangle \\
&- \sum_{j} \Delta t\mathcal H\Big(\{\overline x_{i,j}\},\{\overline p_{i,j}\}, t_j\Big) \Big\}
\end{aligned}
\end{equation} Because we have reversed time in the HJB equation, that change is applied here as well, so the optimal $\{x_{i,j}\}$ in \eqref{eq:HL} is a time-reversed optimal trajectory. Accordingly, we set $\{x_{i,0}\} = \{x_{i,f}\}$ to avoid the infinite penalty assessed by $g(\{x_i\})$, and we set $\{x_{i,J}\} = \{x_i\}$, the point at which we will solve the equation. If we can solve \eqref{eq:HL}, we get an approximation of the value function $u(\{x_i\},t)$ and the optimal trajectories for the agents.

We note the formal similarity between \eqref{eq:HL} and the classical Hopf-Lax formula. Indeed, if $\mathcal H$ is space- and time-independent, then the optimal costate trajectory is constant and the optimal state trajectory is a straight line. In this case, the sums over $j$ telescope, whereupon \eqref{eq:HL} takes the form \eqref{eq:classicHL} once we have stacked $\{x_i\}$ into a single vector variable, so this is indeed a generalization of the classical formula. 

We resolve the saddle-point problem using a version of the primal-dual hybrid gradient method of Chambolle and Pock \cite{PDHG}. Having randomly initialized, the method works by iteratively performing proximal ascent in the costate variable, and then proximal descent and extrapolation in the state variable. Pseudocode for the method is contained in algorithm~\ref{alg:1}. In the algorithm, there is a slight abuse of notation. In the update for $x_{i,j}$, when performing the minimization, we treat $\mathcal H(\{\tilde x_{i}\}, p_{i,j},t_j)$ as a function of the variables $\tilde x_i$ corresponding to the agent $i$ only, keeping variables corresponding to other agents fixed at their values for the current time step $x_{\ell,j}$ for $\ell\neq i$. 

We have taken care to detail algorithm~\ref{alg:1} without combining the state (or costate) variables into a single variable because this demonstrates the manner in which the dependence of $\mathcal H$ on the individual agents is decoupled for several of these computations. For both variables, the minimization regarding variables at different time steps is decoupled. Further, when updating the costate variables, the minimization can be performed at the level of the individual Hamiltonians $H_i$. The situation is less simple for the state variables since the state dynamics are coupled through the collision function $C(\{x_{i}\})$ (and less crucially through the final configuration indicator function $\chi_f(\{x_i\})$). In the next subsection, we give some practical implementation notes for applying the algorithm.

\begin{algorithm}[t!]
\caption{Algorithm for Resolving \eqref{eq:HL}}
Input the point $(\{x_i\},t)$ at which to resolve the HJB equation and proximal step size parameters $\sigma, \tau > 0$ with $\sigma\tau \le 0.25$. 

 Set $\{x_{i,J}\}= \{x_i\}, \{x_{i,0}\}= \{x_{i,f}\}$,  $\{p_{i,0}\} = \{0\}$, and initialize $\{x_{i,j}\}$, $\{p_{i,j}\}$ randomly for all other values. Set $\{z_{i,j}\} = \{x_{i,j}\}$ for all $i,j$. 

\begin{algorithmic}[t!]
\Repeat
 \For {$j = 1$ to $J$ and $i = 1$ to $I$}
		\State \hspace{-2mm}$ \alpha \leftarrow \chi_f(\{x_{i,j}\})C(\{x_{i,j}\})O_i(x_{i,j},t_j)$
    		\State  \hspace{-2mm}$\beta \leftarrow p_{i,j} + \sigma(z_{i,j} - z_{i,j-1})$
    		\State  \hspace{-2mm}$\displaystyle p_{i,j} \leftarrow \argmin_{\tilde p_i \in \R^{n_i}} \{ \Delta t \alpha H_i(t_j,x_{i,j},\tilde p_i) + \tfrac{1}{2\sigma} \lvert \tilde p_i- \beta \rvert^2 \}$
	\EndFor
    
    \State $x^{\text{old}}_{i,j} \leftarrow x_{i,j}$ \,\,\,\,\,\,\,\,\,\,\,\,\,\,\,\,\,\,\, (for each $i,j$)
    \For{$j = 1$ to $J - 1$ and $i = 1$ to $I$ }
    \State \hspace{-2mm}$\xi \leftarrow x_{i,j} - \tau(p_{i,j} - p_{i,j+1})$
    \State \hspace{-2mm}$\displaystyle x_{i,j} \leftarrow \argmin_{\tilde x_i \in \R^{n_i}}\{ - \Delta t \mathcal H(\{\tilde x_i\},\{p_{i,j}\},t_j) + \tfrac{1}{2\tau} \lvert \tilde x_i -\xi \rvert^2 \}$
    \EndFor
    
    \State $z_{i,j} \leftarrow 2x_{i,j} -x^{\text{old}}_{i,j}$  \,\,\,\,\, (for each $i,j$)
    
\Until{convergence}

\State $u =\sum_{i,j} \langle p_{i,j}, x_{i,j} - x_{i,j-1} \rangle - \Delta t \sum_j \mathcal H(\{x_{i,j}\},\{p_{i,j}\},t_j)$
\State \textbf{return } $u$, $\{x_{i,j}\}$
\end{algorithmic}
\label{alg:1}
\end{algorithm}

\subsection{Practical Implementation Notes}

As a first note, we mention that for many nontrivial models of motion, one can resolve the optimization in the individual Hamiltonians $H_i$ defined in \eqref{eq:Hi} so that the Hamiltonian is given by an explicit formula. Two common examples are that of isotropic motion or simple vehicles. For an agent exhibiting isotropic motion in $\mathbb R^n$, the dynamics are $\dot {\B x} = V(\B x)\B a$, where $V(x)>0$ is a local speed function and the control variable $\B a(\cdot) \in \mathbb S^{n-1}$ represents the direction of motion. In this case, an easy computation shows that $H(x,p) = V(x)\abs{p}$. Likewise, a common model for a simple car considers $X = (x,y,\theta)$ satisfying $\dot x = v\cos(\theta), \dot y = v\sin(\theta), \dot \theta = \omega$. Here $(x,y)$ denote the center of mass of the car and $\theta$ is its orientation, and the control variables $v \in [-V,V],\omega \in [-W,W]$ correspond to tangential and angular velocity respectively (with $V,W$ being bounds on these quantities). In this case, the Hamiltonian is $H(X,p) = V\abs{\cos(\theta)p_1 + \sin(\theta)p_2} + W\abs{p_3}$ (this Hamiltonian and higher dimensional versions are derived in \cite{TakeiTsai2}). For a slightly more exotic example, we also consider a quadcopter drone whose dynamics are \begin{equation}
\label{eq:droneDyn} \begin{aligned}
&\ddot x=v(\sin\phi\sin\psi+\cos\phi\cos\psi\sin\theta)\\
    &\ddot y=v(\cos\phi\sin\theta\sin\psi-\cos\psi\sin\phi)\\
    &\ddot z=v\cos\theta\cos\phi-g\\
    &\ddot\psi=\tau_\psi, \,\,\,\,\,\,
    \ddot\theta=\tau_\theta, \,\,\,\,\,\,
    \ddot\phi=\tau_\phi.
\end{aligned}
\end{equation} The full derivation of these equations is in \cite{Carrillo}. For our purposes, it is enough to understand that $(x,y,z)$ represent the center of mass of the drone and $(\psi,\theta,\phi)$ its orientation. The control variables here are the upward thrust $v$ and angular thrusts $\tau_\psi,\tau_\theta,\tau_\phi$. Each of these is normalized to take values in $[-1,1]$ and the mass of the drone is normalized to $1$. The parameter $g$ in \eqref{eq:droneDyn} is the downward accelation due to gravity. We note that in this case, the control takes place at the level of acceleration rather than velocity. Thus setting $X = (x,y,z,\psi,\theta,\phi)$, the dynamics are 12-dimensional with state space comprised  of variables $(X,\dot X)$. Writing the costate variable similarly as $(P,\dot P)$, the Hamiltonian is \begin{equation}\label{eq:DroneH} \small
H(X,\dot X,P,\dot P) = -\langle \dot X, P\rangle +\abs{\langle \dot P_{1:3},\gamma(X_{4:6})\rangle} + g\dot P_3+\abs{\dot P_{4:6}}_1
\end{equation}  where  $\gamma(X_{4:6}) = \gamma(\psi,\theta,\phi)$ is the unit vector defined from \eqref{eq:droneDyn} so that $(\ddot x,\ddot y,\ddot z) = v\gamma(\psi,\theta,\phi) - (0,0,g)$. 

In any case, because one needs to resolve roughly $IJ$ optimization problems in each iteration of algorithm~\ref{alg:1}, in order to maintain efficiency it is important to resolve these exactly using analytic formulas whenever possible. For the Hamiltonians listed above, one can explicitly write down the solutions of the optimization problems involving costate variables $\{p_{i,j}\}$. In the case of isotropic motions, the relevant formula is included in \cite{ParkinsonPolage} and in the case of simple vehicles in \cite{ParkinsonBoyle}. For the drone dynamics with Hamiltonian \eqref{eq:DroneH}, using the state variables $(X,\dot X)$ and costate variables $(P,\dot P)$ (and writing the auxiliary vector $\beta$ in the first loop of algorithm~\ref{alg:1} analogously as $(\beta,\dot \beta)$), the explicit update rule for the costate vectors is \begin{equation}\label{eq:DroneP} \begin{aligned}P &\leftarrow \alpha \dot X + \beta, \\ 
\dot P_{1:3} &\leftarrow -\min\Big(1,\tfrac{\alpha \sigma \Delta t}{\abs{\langle \gamma,B\rangle}}\Big)\langle \gamma,B\rangle\gamma + B, \\
\dot P_{k} &\leftarrow \max\Big(0,1-\tfrac{\alpha \sigma \Delta t}{\abs{\dot \beta_k}}\Big)\dot \beta_k, \,\,\,\,\,\,\,\,\, k =4,5,6,
\end{aligned}\end{equation} where $\alpha$ is as in algorithm~\ref{alg:1}, $\gamma$ is as in \eqref{eq:DroneH}, and $B = \dot \beta_{1:3} - \alpha \sigma \Delta t(0,0,g)$. This update is applied to each $p_{i,j}$ in the first loop in algorithm~\ref{alg:1}. 

The optimization problems which update the vectors $x_{i,j}$ cannot typically be resolved analytically. However, it is observed in \cite{Lin} (and corroborated in \cite{ParkinsonBoyle,ParkinsonPolage,ParkinsonBaca}) that these can approximated quite crudely at each iteration and algorithm~\ref{alg:1} will still produce a good approximation of the optimal trajectories. Thus, following the suggestion of \cite{Lin}, we simply approximate the minimizer which defines $x_{i,j}$ in algorithm~\ref{alg:1} using a single step of gradient descent, while using smooth approximations of any indicator functions. Specifically, we use the approximations \begin{equation} \label{eq:As}\begin{aligned} \chi_f(\{x_i\}) &= 1-e^{-A_1\sum_i \abs{x_i-x_{i,f}}^2},\\
c_{k\ell}(x_k,x_\ell) &= \tfrac 1 2 \Big(1+\tanh\big(A_2(\abs{x_k-x_\ell}^2 -\delta^2)\big)\Big),\\
O_i(x_i,t) &= \tfrac 1 2 \Big(1+\tanh\big(A_3d(x_i,t)^2\big)\Big),
\end{aligned}\end{equation} for some large constants $A_1,A_2,A_3$, where $\delta > 0$ is the collision radius and $d(x,t)$ is the signed distance function to the obstacles at time $t$ (negative inside the obstacles). Thus $c_{k\ell}(x_k,x_\ell) \approx 0$ if agents $k$ and $\ell$ are within $\delta$ of each other, and $c_{k\ell}(x_k,x_\ell) \approx 1$ if agents $k$ and $\ell$ are sufficiently far apart. Likewise, $O_i(x_i,t) \approx 0$ if $x_i$ is in an obstacle and $O_i(x_i,t)\approx 1$ if $x_i$ is in free space. We note that there is another slight abuse of notation here: $c_{k\ell}$ and $O_i$ operate only on the \emph{spatial} components of the dynamics. For example, these functions will ignore the velocities $\dot x,\dot y,\dot z$ and the angular components for the drone dynamics \eqref{eq:droneDyn}. Finally, for quick computation of the obstacle functions, we use the strategy of \cite{ParkinsonBoyle} wherein obstacles are required to have simple geometry, but if desired, more complicated geometries can be approximated by, for example, a collection of disjoint balls. 

\section{Results \& Discussion}

In this section, we demonstrate the efficacy of our model and numerical scheme with some examples. All simulations were run on the first author's desktop computer which has an Intel(R) Core(TM) i7-14700 processor and 32GB of RAM. In all cases, we set $\Delta t = 0.1$, and the parameters from algorithm~\ref{alg:1} to $\sigma = 1, \tau = 0.25$. For the large constants in \eqref{eq:As}, we set $A_2=A_3 = 100$, and $A_1 = 10$ at the outset. We then set $A_1\leftarrow \min(A_1 + 50,1000)$ after each 1000 iterations. Likewise, we approximate the minimization for $\{x_{i,j}\}$ with a single step of gradient descent with an initial rate of $0.1$ and we halve this rate after each $1000$ iterations as the path refines. Our examples focus on simple cars like those in \cite{ParkinsonBoyle} with maximum tangential velocity $1$ and maximum angular velocity $2$, as well as quadcopter drones with maximum upward acceleration normalized to $1$ and downward acceleration due to gravity set at $0.1$, which is consistent which high performance drones which have upward acceleration exceeding 10 gs \cite{Bond}. The code for these figures (which also produces animated versions) is available at the link below.\footnote{\href{https://github.com/chparkinson/Collisionless-MultiAgent-HJB}{\texttt{github.com/chparkinson/Collisionless-MultiAgent-HJB}}}

\begin{figure}[b!]
\centering
{\setlength{\fboxrule}{1.5pt}
\fbox{\includegraphics[width=0.48\textwidth]{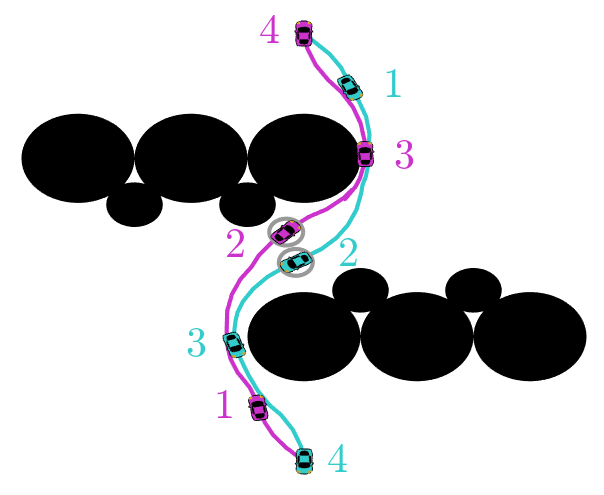}}}
\caption{Two cars navigate around obstacles to their final destinations. Notice they are careful to narrowly avoid colliding at time $2$.}
\label{fig:1}
\end{figure}

Our first example in figure~\ref{fig:1} is a rather simplistic proof of concept. We have two cars navigating a domain with several obstacles (black). The first car (cyan) begins at the top of the frame and must end at the bottom, and vice-versa for the second (magenta). We have plotted their paths in a single image, labeling successive points along the path. Specifically, at the time labeled $2$, we have drawn a circle around each car of radius $\delta/2$, meaning that if these circles overlap, the cars are considered to have collided. However, as seen, the cars pass each other at exactly the allowable distance to avoid a collision. In this example, the state space is 6 dimensional and the computation requires on the order of 2 seconds of CPU time and 2000 iterations in algorithm~\ref{alg:1}.

\begin{figure*}[!]
\centering
\fbox{\includegraphics[width=0.18\textwidth,trim = 70 25 50 20, clip]{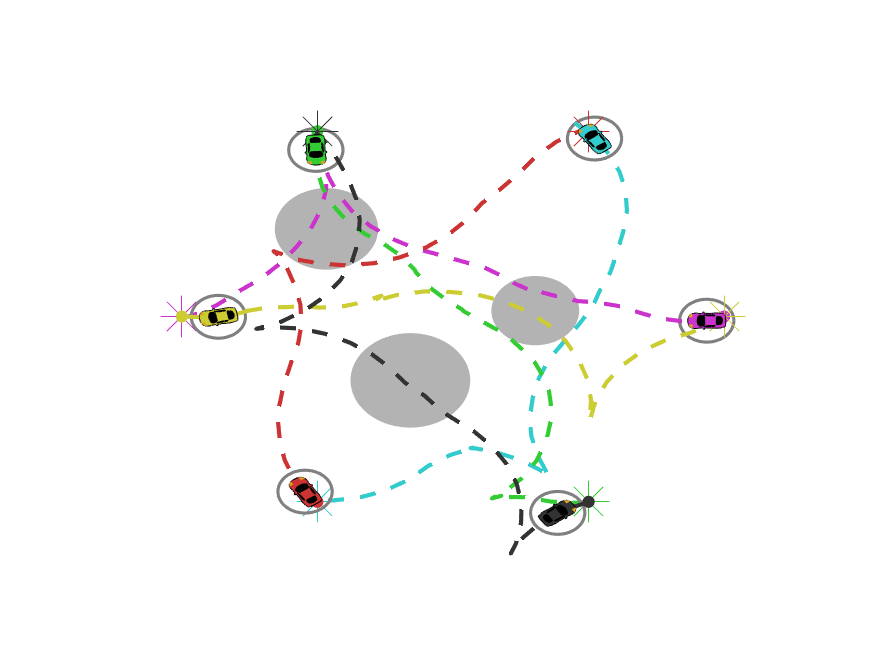}} \, 
\fbox{\includegraphics[width=0.18\textwidth,trim = 70 25 50 20, clip]{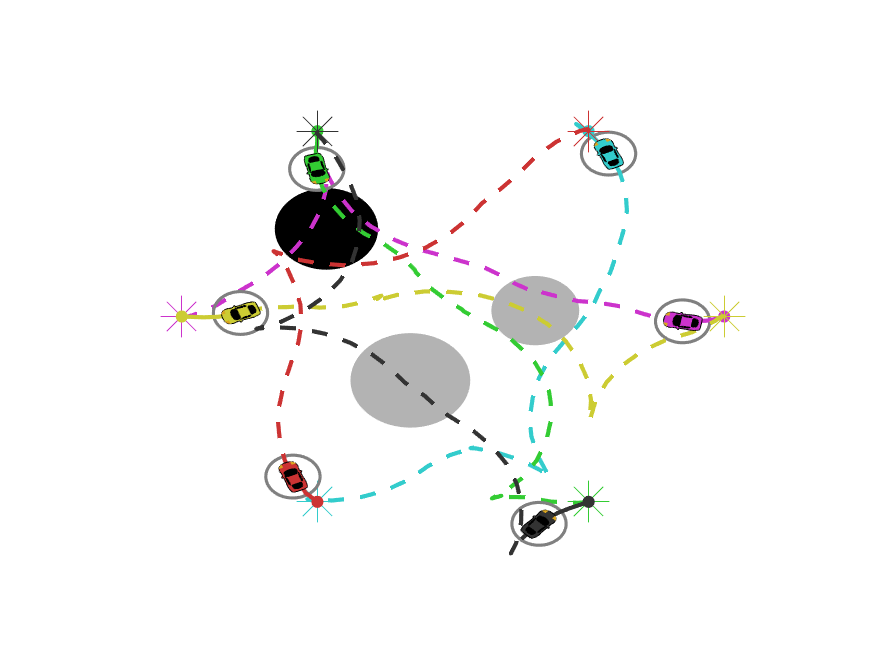}} \, 
\fbox{\includegraphics[width=0.18\textwidth,trim = 70 25 50 20, clip]{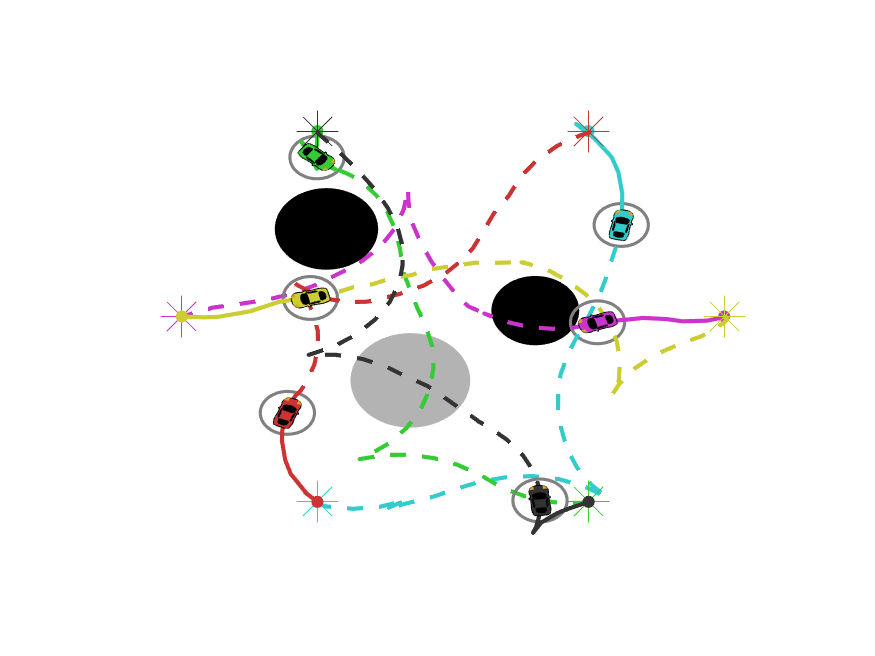}}\, 
\fbox{\includegraphics[width=0.18\textwidth,trim = 70 25 50 20, clip]{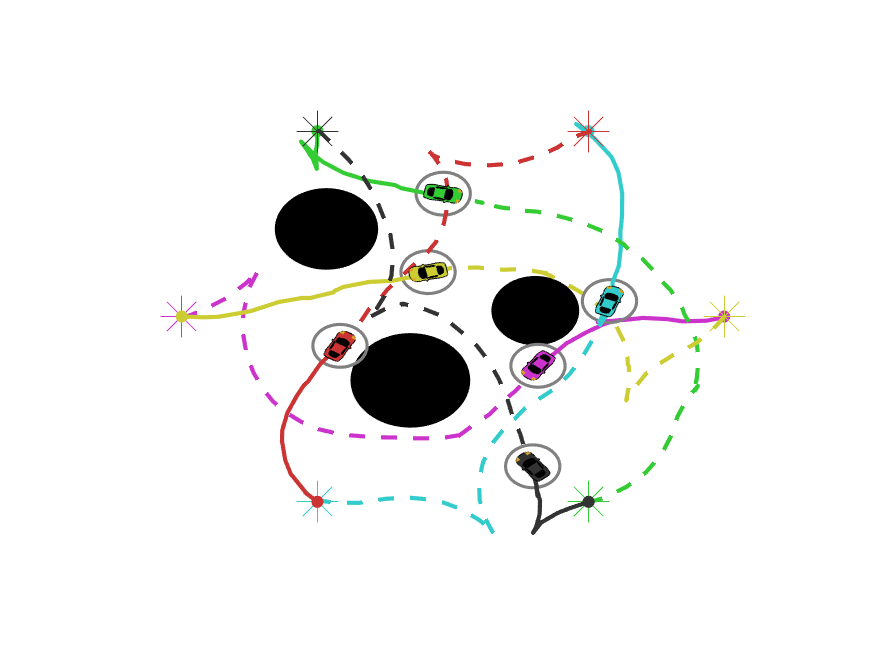}} \, 
\fbox{\includegraphics[width=0.18\textwidth,trim = 70 25 50 20, clip]{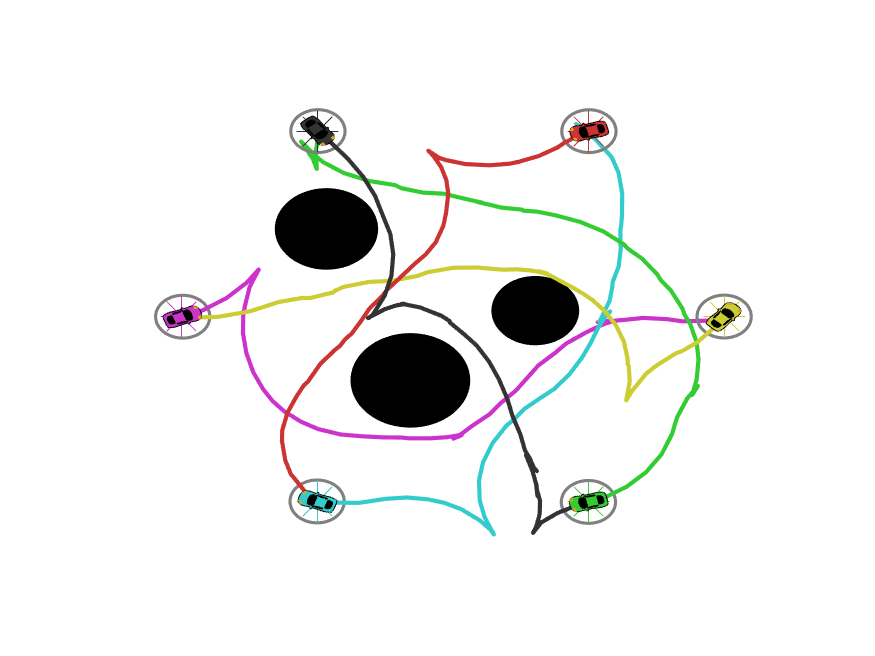}} 
\caption{Six cars navigate to antipodal points around a circle. Gray circles are undiscovered obstacles, black circles are discovered obstacles. Dotted lines are the currently resolved paths. These are recomputed to avoid obstacles as needed. }
\label{fig:2}
\end{figure*}

Our second example---included in figure~\ref{fig:2}---integrates the method with a semi-real-time corrector when new information arises. In particular, we have six cars evenly spaced around a circle, each needing to travel to the antipodal point while avoiding each other and obstacles. At the outset, the location of obstacles is unknown. When a car discovers an obstacle, we recompute the optimal paths with this new information. In the series of pictures, gray circles are undicovered obstacles, black circles are discovered obstacles, and we have paused at the moments when the obstacles are discovered. The dotted paths are the current paths the cars are planning to take, but as we see, these sometimes intersect with unseen obstacles and are eventually recomputed. In this case, the state space is 18 dimensional and there are 3 obstacles all of are discovered, requiring paths to be recomputed, so a total of 4 optimal paths are resolved. Each path required roughly 3.5 seconds of CPU time and fewer than 1000 iterations to resolve, so the total CPU time for this example was roughly 14 seconds. This demonstrates that, even for relatively high dimensional state space, our algorithm is efficient enough to be implemented in semi-real-time. We note that the plotting of the quadcopter drone was adapted from \cite{Singh}.

\begin{figure*}[!]
\centering
\fbox{\includegraphics[width=0.18\textwidth,trim = 70 25 50 20, clip]{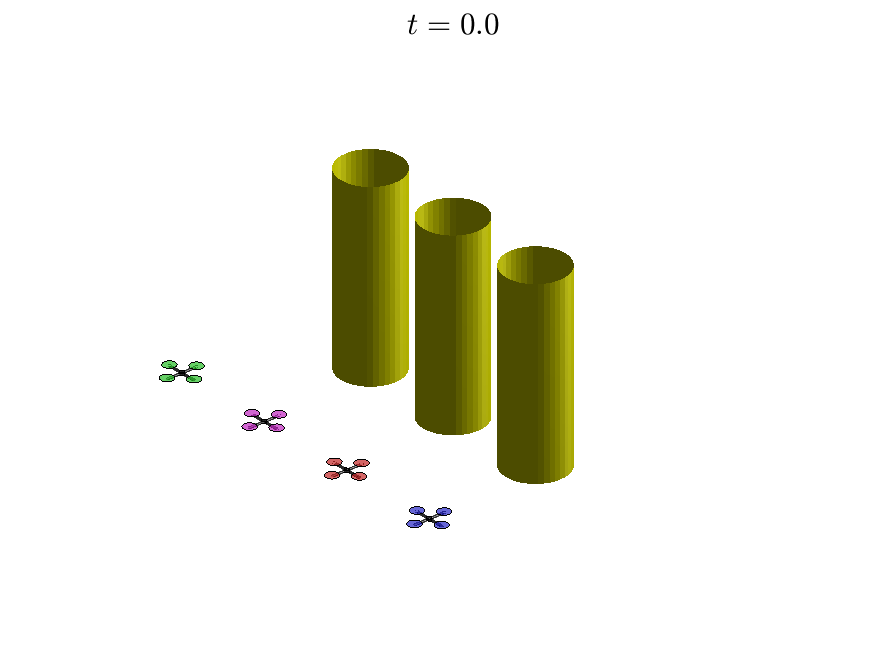}} \, 
\fbox{\includegraphics[width=0.18\textwidth,trim = 70 25 50 20, clip]{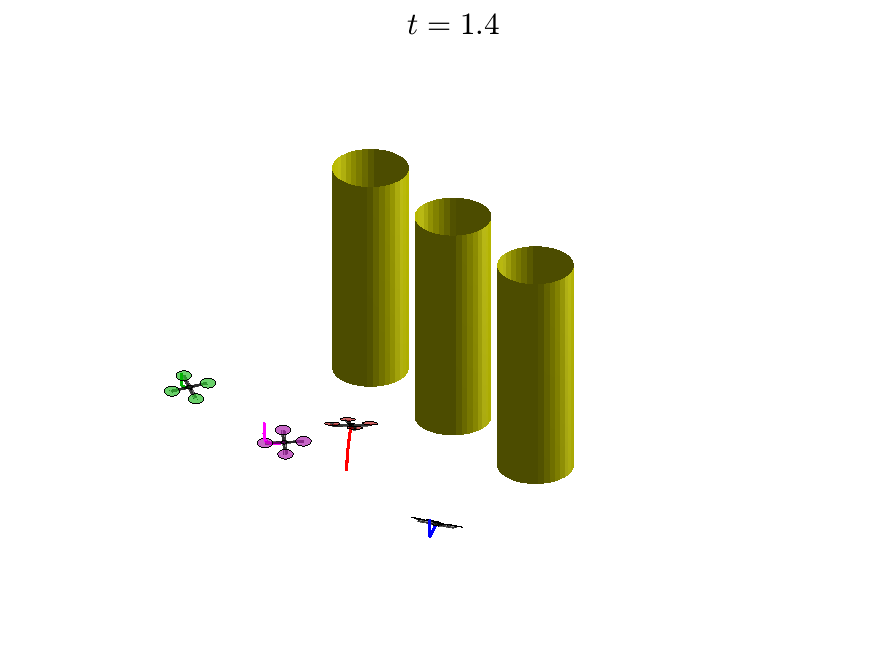}} \, 
\fbox{\includegraphics[width=0.18\textwidth,trim = 70 25 50 20, clip]{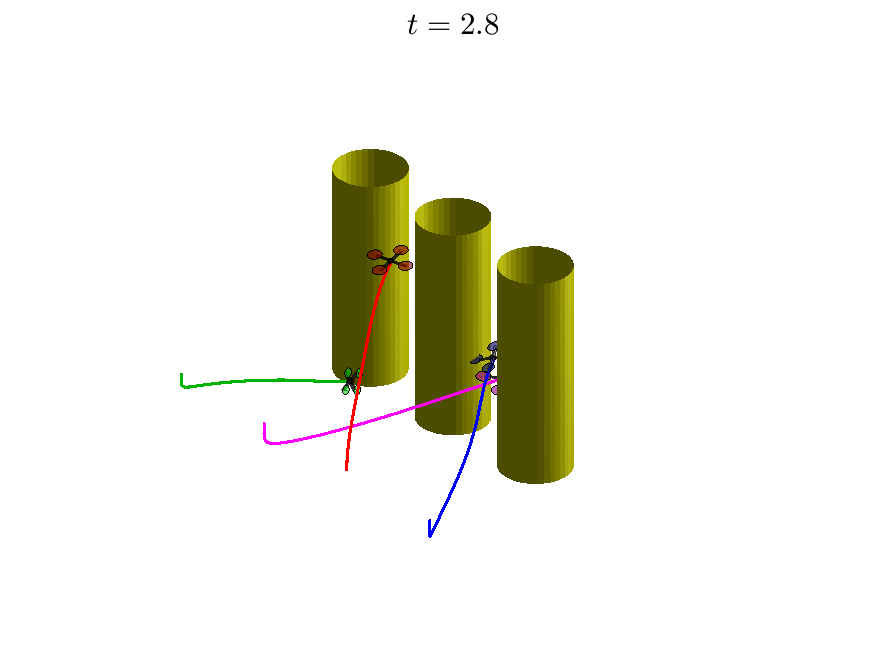}}\, 
\fbox{\includegraphics[width=0.18\textwidth,trim = 70 25 50 20, clip]{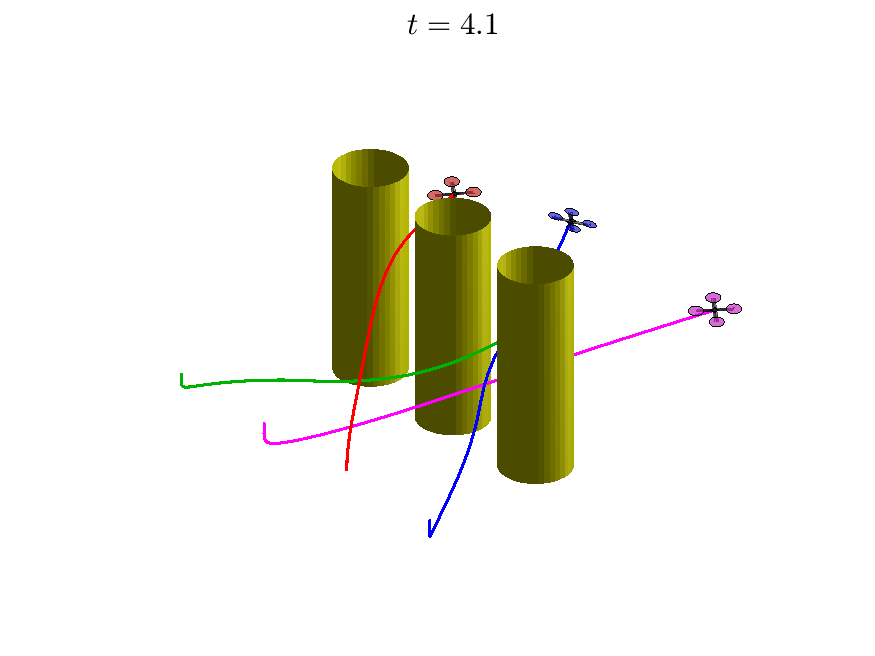}} \, 
\fbox{\includegraphics[width=0.18\textwidth,trim = 70 25 50 20, clip]{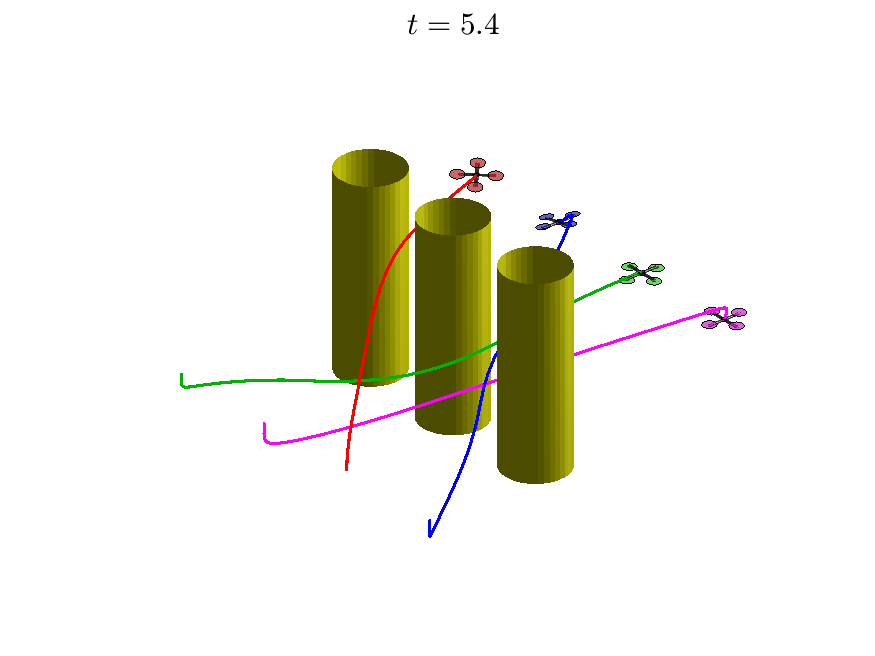}} 
\caption{Four quadcopter drones navigate around obstacles and each other to their final destinations. Here the obstacles are the yellow cylinders. Notice that the drones on the left in the initial picture (green and pink lines) allow themselves to drop a bit at the beginning, so that they can then fly unobstructed between the obstacles. This is to avoid collisions which are harder to prevent in a dynamic model compared to a kinematic model, and hence need to be considered long in advance.}
\label{fig:3}
\end{figure*}

Our last example involves four quadcopter drones navigating around obstacles. This is included in figure~\ref{fig:3}. Here the obstacles are the yellow cylinders. We particularly note that the green and pink drones (left in initial picture) do not start immediately, instead opting to let themselves drop a bit. This is to avoid collisions when going through the narrow corridor between obstacles. Because this a dynamical model, wherein control occurs at the level of acceleration, it is important to plan for collision avoidance long in advance, since local adjustments to paths are more difficult to accomplish. The state space in this example is 48 dimensions and algorithm~\ref{alg:1} typically requires roughly 10 seconds of CPU time and 3000 iterations to resolve optimal paths. This demonstrates that our model can accommodate reasonably complex dynamics. 

\bibliographystyle{ieeetr}
\bibliography{Biblio}

\begin{thebibliography}{10}

\bibitem{Evac}
L.~Sun, H.~Liu, B.~Fan, X.~Li, and W.~Li, ``Multi-agent deep reinforcement
  learning-based path planning for crowd evacuation in buildings,'' {\em
  Journal of Building Engineering}, vol.~111, p.~113345, 2025.

\bibitem{AirTravel}
A.~Ronzhin, A.~Saveliev, D.~Anikin, A.~Zaytseva, E.~Cherskikh, and A.~Figurek,
  ``Collision-free multi-uav 3d path planning using dynamic spatial
  reservations of reducing air corridors,'' {\em Robotics and Autonomous
  Systems}, p.~105359, 2026.

\bibitem{Coverage}
X.-W. Ma, T.~Huang, W.-L. Liu, and Y.-J. Gong, ``Collision-aware evolutionary
  algorithm for multi-agent coverage path planning,'' in {\em 2024 11th
  International Conference on Machine Intelligence Theory and Applications
  (MiTA)}, pp.~1--8, IEEE, 2024.

\bibitem{DL1}
R.~Alharthi, I.~Noreen, A.~Khan, T.~Aljrees, Z.~Riaz, and N.~Innab, ``Novel
  deep reinforcement learning based collision avoidance approach for path
  planning of robots in unknown environment,'' {\em PloS one}, vol.~20, no.~1,
  p.~e0312559, 2025.

\bibitem{DL2}
T.~T. Nguyen, S.~Nahavandi, I.~Razzak, D.~Nguyen, N.~T. Pham, and Q.~V.~H.
  Nguyen, ``The emergence of deep reinforcement learning for path planning,''
  in {\em 2025 IEEE International Conference on Systems, Man, and Cybernetics
  (SMC)}, pp.~6265--6272, IEEE, 2025.

\bibitem{DL3}
Z.~Wang, S.~Song, and S.~Cheng, ``Path planning of mobile robot based on
  improved double deep q-network algorithm,'' {\em Frontiers in Neurorobotics},
  vol.~19, p.~1512953, 2025.

\bibitem{Robust1}
C.~Eleftheriadis, A.~Symeonidis, and P.~Katsaros, ``Adversarial robustness
  improvement for deep neural networks,'' {\em Machine Vision and
  Applications}, vol.~35, no.~3, p.~35, 2024.

\bibitem{Explain1}
A.~Ellouze, M.~Karray, and M.~Ksantini, ``A hybrid decision-making framework
  for autonomous vehicles in urban environments based on multi-agent
  reinforcement learning with explainable ai,'' {\em Vehicles}, vol.~8, no.~1,
  p.~8, 2026.

\bibitem{Robust2}
J.~Li and G.~Li, ``Triangular trade-off between robustness, accuracy, and
  fairness in deep neural networks: A survey,'' {\em ACM Computing Surveys},
  vol.~57, no.~6, pp.~1--40, 2025.

\bibitem{Explain2}
G.~Zhang, J.~Song, and J.~Qiao, ``Explainable ai in deep learning: Methods,
  challenges, and future directions,'' {\em Authorea Preprints}, 2026.

\bibitem{Lolla}
T.~Lolla, M.~P. Ueckermann, K.~Yi{\u{g}}it, P.~J. Haley, and P.~F. Lermusiaux,
  ``Path planning in time dependent flow fields using level set methods,'' in
  {\em 2012 IEEE International Conference on Robotics and Automation},
  pp.~166--173, IEEE, 2012.

\bibitem{Arnold}
D.~J. Arnold, D.~Fernandez, R.~Jia, C.~Parkinson, D.~Tonne, Y.~Yaniv, A.~L.
  Bertozzi, and S.~J. Osher, ``Modeling environmental crime in protected areas
  using the level set method,'' {\em SIAM Journal on Applied Mathematics},
  vol.~79, no.~3, pp.~802--821, 2019.

\bibitem{Cartee2}
E.~Cartee and A.~Vladimirsky, ``Control-theoretic models of environmental
  crime,'' {\em SIAM Journal on Applied Mathematics}, vol.~80, no.~3,
  pp.~1441--1466, 2020.

\bibitem{Chen}
B.~Chen, K.~Peng, C.~Parkinson, A.~L. Bertozzi, T.~L. Slough, and
  J.~Urpelainen, ``Modeling illegal logging in {Brazil},'' {\em Research in the
  Mathematical Sciences}, vol.~8, no.~2, pp.~1--21, 2021.

\bibitem{Parkinson}
C.~Parkinson, D.~Arnold, A.~L. Bertozzi, Y.~T. Chow, and S.~Osher, ``Optimal
  human navigation in steep terrain: a {{H}amilton-{J}acobi-Bellman}
  approach,'' {\em Communications in Mathematical Sciences}, vol.~17, no.~1,
  pp.~227--242, 2019.

\bibitem{Parkinson2}
C.~Parkinson, D.~Arnold, A.~Bertozzi, and S.~Osher, ``A model for optimal human
  navigation with stochastic effects,'' {\em SIAM Journal on Applied
  Mathematics}, vol.~80, no.~4, pp.~1862--1881, 2020.

\bibitem{Cartee1}
E.~Cartee, L.~Lai, Q.~Song, and A.~Vladimirsky, ``Time-dependent
  surveillance-evasion games,'' in {\em 2019 IEEE 58th Conference on Decision
  and Control (CDC)}, pp.~7128--7133, IEEE, 2019.

\bibitem{GeeFear}
M.~Gee, N.~Gonzalez-Granda, S.~Joshi, N.~Rudrapatna, A.~Somalwar, S.~P. Ellner,
  and A.~Vladimirsky, ``Navigating the landscape of fear,'' {\em bioRxiv},
  pp.~2024--08, 2024.

\bibitem{TakeiTsai1}
R.~Takei, R.~Tsai, H.~Shen, and Y.~Landa, ``A practical path-planning algorithm
  for a simple car: a {H}amilton-{J}acobi approach,'' in {\em Proceedings of
  the 2010 American Control Conference}, pp.~6175--6180, June 2010.

\bibitem{TakeiTsai2}
R.~Takei and R.~Tsai, ``Optimal trajectories of curvature constrained motion in
  the {H}amilton-{J}acobi formulation,'' {\em Journal of Scientific Computing},
  vol.~54, pp.~622--644, Feb 2013.

\bibitem{ParkinsonCar1}
C.~Parkinson, A.~L. Bertozzi, and S.~J. Osher, ``A {{H}amilton-{J}acobi}
  formulation for time-optimal paths of rectangular nonholonomic vehicles,'' in
  {\em 2020 59th IEEE Conference on Decision and Control (CDC)},
  pp.~4073--4078, IEEE, 2020.

\bibitem{ParkinsonCar2}
C.~Parkinson and M.~Ceccia, ``Time-optimal paths for simple cars with moving
  obstacles in the {H}amilton-{J}acobi formulation,'' in {\em 2022 American
  Control Conference (ACC)}, pp.~2944--2949, IEEE, 2022.

\bibitem{Gee}
M.~Gee and A.~Vladimirsky, ``Optimal path-planning with random breakdowns,''
  {\em IEEE Control Systems Letters}, vol.~6, pp.~1658--1663, 2021.

\bibitem{FS1}
S.~Luo and H.~Zhao, ``Convergence analysis of the fast sweeping method for
  static convex {H}amilton--{J}acobi equations,'' {\em Research in the
  Mathematical Sciences}, vol.~3, no.~1, p.~35, 2016.

\bibitem{FS4}
C.~Parkinson, ``A rotating-grid upwind fast sweeping scheme for a class of
  {H}amilton-{J}acobi equations,'' {\em Journal of Scientific Computing},
  vol.~88, no.~1, p.~13, 2021.

\bibitem{Sethian}
J.~A. Sethian, ``A fast marching level set method for monotonically advancing
  fronts,'' {\em Proceedings of the National Academy of Sciences}, vol.~93,
  no.~4, pp.~1591--1595, 1996.

\bibitem{SethVlad1}
J.~A. Sethian and A.~Vladimirsky, ``Ordered upwind methods for static
  {H}amilton-{J}acobi equations: Theory and algorithms,'' {\em SIAM Journal on
  Numerical Analysis}, vol.~41, no.~1, pp.~325--363, 2003.

\bibitem{HighDimHJ}
T.~Meng, S.~Liu, S.~W. Fung, and S.~Osher, ``Recent advances in numerical
  solutions for hamilton-jacobi pdes,'' {\em arXiv preprint arXiv:2502.20833},
  2025.

\bibitem{Lin}
A.~T. Lin, Y.~T. Chow, and S.~J. Osher, ``A splitting method for overcoming the
  curse of dimensionality in {H}amilton--{J}acobi equations arising from
  nonlinear optimal control and differential games with applications to
  trajectory generation,'' {\em Communications in Mathematical Sciences},
  vol.~16, 1 2018.

\bibitem{ParkinsonBoyle}
C.~Parkinson and I.~Boyle, ``Efficient and scalable path-planning algorithms
  for curvature constrained motion in the {H}amilton-{J}acobi formulation,''
  {\em Journal of Computational Physics}, vol.~509 (113050), 2024.

\bibitem{ParkinsonPolage}
C.~Parkinson and K.~Polage, ``An efficient semi-real-time algorithm for path
  planning in the {H}amilton--{J}acobi formulation,'' {\em IEEE Control Systems
  Letters}, vol.~7, pp.~3621--3626, 2023.

\bibitem{ParkinsonBaca}
C.~Parkinson and A.~Baca, ``A hopf-lax type formula for multi-agent path
  planning with pattern coordination,'' in {\em 2025 IEEE 64th Conference on
  Decision and Control (CDC)}, pp.~649--654, IEEE, 2025.

\bibitem{Bardi}
M.~Bardi, I.~C. Dolcetta, {\em et~al.}, {\em Optimal control and viscosity
  solutions of Hamilton-Jacobi-Bellman equations}, vol.~12.
\newblock Springer, 1997.

\bibitem{Tran}
H.~V. Tran, {\em Hamilton--Jacobi equations: theory and applications},
  vol.~213.
\newblock American Mathematical Soc., 2021.

\bibitem{STLC}
M.~I. Krastanov and M.~N. Nikolova, ``On the small-time local
  controllability,'' {\em Systems \& Control Letters}, vol.~177, p.~105535,
  2023.

\bibitem{PDHG}
A.~Chambolle and T.~Pock, ``A first-order primal-dual algorithm for convex
  problems with applications to imaging,'' {\em Journal of mathematical imaging
  and vision}, vol.~40, no.~1, pp.~120--145, 2011.

\bibitem{Carrillo}
L.~R.~G. Carrillo, A.~E.~D. Lopez, R.~Lozano, and C.~Pegard, {\em Quad
  Rotorcraft Control}.
\newblock Springer London, Mar 2018.

\bibitem{Bond}
E.~Bond, B.~Crowther, and B.~Parslew, ``The rise of high-performance
  multi-rotor unmanned aerial vehicles-how worried should we be?,'' in {\em
  2019 Workshop on Research, Education and Development of Unmanned Aerial
  Systems (RED UAS)}, pp.~177--184, IEEE, 2019.

\bibitem{Singh}
J.~Singh, ``Quadcopter model: {MATLAB} code for animation.'' MATLAB Central
  File Exchange.
  \url{https://www.mathworks.com/matlabcentral/fileexchange/97192-quadcopter-model-matlab-code-for-animation},
  2024.

\end{thebibliography}

\end{document}